\documentclass[draft,a4paper]{amsart}
\title{ Discrete Logarithms in Generalized Jacobians }
\author{ S.~D.~Galbraith \and B.~A.~Smith }
\address{Mathematics Department, Royal Holloway University
of London, Egham, Surrey TW20 0EX, UK.}
\email{Steven.Galbraith@rhul.ac.uk, Ben.Smith@rhul.ac.uk}

\usepackage{a4wide}
\usepackage{amsmath,amsfonts,amssymb,amsthm}

\theoremstyle{definition}
\newtheorem{definition}{Definition}[section]

\theoremstyle{plain}

\newtheorem{proposition}[definition]{Proposition}
\newtheorem{lemma}[definition]{Lemma}
\newtheorem{corollary}[definition]{Corollary}
\theoremstyle{remark}
\newtheorem{example}[definition]{Example}
\newtheorem{remark}[definition]{Remark}

\newcommand{\multiplicativegroup}{\mathbb{G}_{m}}

\newcommand{\modulus}{\mathfrak{m}}

\newcommand{\grp}[1]{{\mathcal{#1}}}
\newcommand{\Jac}[1]{{\ensuremath{\grp{J}_{#1}}}}

\begin{document}

\begin{abstract}
D\'ech\`ene has proposed generalized Jacobians 
as a source of groups for public-key cryptosystems
based on the hardness of the Discrete Logarithm Problem (DLP).
Her specific proposal gives rise to a group isomorphic to the
semidirect product of an elliptic curve and a multiplicative
group of a finite field.
We explain why her proposal has no advantages 
over simply taking the direct product of groups.
We then argue that generalized Jacobians offer poorer security 
and efficiency than standard Jacobians.
\end{abstract}

\maketitle

\section{Introduction}

Recently,
D\'ech\`ene~\cite{Dechene-ANTS}
has proposed generalized Jacobians 
as a source of groups for public-key cryptosystems
based on the hardness of the Discrete Logarithm Problem (DLP).
Generalized Jacobians offer a natural generalization of both 
torus-based and curve-based cryptography.

D\'ech\`ene's specific proposal gives rise to a group isomorphic to the
semidirect product of an elliptic curve $\grp{E}(k)$
and a multiplicative group of a finite field $\multiplicativegroup(k)$.
She remarks in Section 6 of~\cite{Dechene-ANTS}
that the DLP in such a generalized Jacobian
can be reduced to sequentially solving a DLP in $\grp{E}(k)$
followed by a DLP in $\multiplicativegroup(k)$
and so it is ``at least as hard as a DLP in $\grp{E}(k)$
and at least as hard as a DLP in $\multiplicativegroup(k)$''.

Our main observation follows from applying the standard
Pohlig-Hellman reduction and therefore
reducing to the case of elements of prime order.
It then immediately follows (see Proposition~\ref{DLP-reduction})
that one can solve the DLP in the generalized
Jacobian by solving a number of DLPs in 
$\grp{E}(k)$ and $\multiplicativegroup(k)$ in parallel.
One concludes that the generalized Jacobian DLP is
\textbf{at most} as hard as the DLP in $\grp{E}(k)$ and 
the DLP in $\multiplicativegroup(k)$.
As we will explain, one can get the same
security with greater efficiency by simply taking
the direct product $\grp{E}(k) \times \multiplicativegroup(k)$.


In our presentation
we consider the DLP in the simpler and more general setting
of extensions of algebraic groups.
We will argue that extensions offer no advantages 
over the existing Jacobian or torus constructions
for DLP-based cryptography.

\bigskip

Throughout this article,
we let $k$ be a finite field.
All varieties are nonsingular $k$-varieties.
We say that a morphism of algebraic groups is {\em explicit}
if it may be evaluated in polynomial time.
Algebraic groups are said to be {\em explicitly isomorphic}
if there is an explicit isomorphism between them.
All algebraic groups in this article are commutative,
and written additively.
We denote algebraic groups with script letters
and their underlying varieties with capital letters:
so if $\grp{A}$ is an algebraic group,
then $A$ denotes its underlying variety.

\section{Discrete Logarithms in Extensions of Commutative Algebraic Groups}

Fix a pair of algebraic groups $\grp{A}$ and $\grp{B}$.
An {\em extension} of $\grp{A}$ by $\grp{B}$ is an algebraic group $\grp{C}$
together with separable homomorphisms $\iota: \grp{B} \to \grp{C}$ and $\pi: \grp{C} \to \grp{A}$,
all defined over $k$,
such that the following sequence is exact:
\begin{equation}
\label{GE-SES}
	0 \to 
	\grp{B} \stackrel{\iota}{\longrightarrow}
	\grp{C} \stackrel{\pi}{\longrightarrow}
	\grp{A} \to 
	0 .
\end{equation}
We will assume that the maps $\iota$, $\pi$, and $\iota^{-1}$ 
(where it is defined)
are explicit.
A trivial example of an extension of $\grp{A}$ by $\grp{B}$
is the direct product $\grp{C} = \grp{A}\times \grp{B}$,
with $\iota$ and $\pi$ the obvious maps.
The motivating example for this work
is the case where $\grp{C}$ is a generalized Jacobian:
here $\grp{A}$ is the Jacobian of an algebraic curve,
$\grp{B}$ is a certain affine algebraic group,\footnote{
	The algebraic group in question
	is isomorphic to a product of multiplicative groups
	(i.e.~a torus),
	together with a product of Witt groups,
	in which the DLP is trivial.
}
and the group structure of $\grp{C}$ is determined by a map
$c_\modulus: \grp{A}^2 \to \grp{B}$.
The generalized Jacobians proposed for cryptography by D\'ech\`ene
are the special case where $\grp{A}$ is an elliptic curve
and $\grp{B}$ is the multiplicative group.

We wish to assess the suitability of $\grp{C}$
as a source of groups for cryptography,
compared with $\grp{A}$ and $\grp{B}$.
Suppose we wish to solve a DLP in a subgroup $\grp{G}$ of $\grp{C}(k)$.
The group $\grp{G}$ is necessarily finite,
and without loss of generality
we may assume that $\grp{G}$ is cyclic.
By the standard reduction of Pohlig and Hellman \cite{Pohlig-Hellman},
we may reduce to the case where the order of $\grp{G}$ is prime.

\begin{proposition}
\label{DLP-reduction}
	Let $\grp{G}$ be a subgroup of $\grp{C}(k)$, of prime order $l$.
	If $\grp{G}$ is contained in $\iota(\grp{B})$,
	then the DLP in $\grp{G}$ reduces to 
	the DLP in a subgroup of order $l$ in $\grp{B}(k)$.
	Otherwise,
	the DLP in $\grp{G}$ reduces to 
	the DLP in a subgroup of order $l$ in $\grp{A}(k)$.
\end{proposition}
\begin{proof}
	If $\grp{G}$ is a subgroup of $\iota(\grp{B})$,
	then it is explicitly isomorphic to 
	the subgroup $\iota^{-1}(\grp{G})$ of $\grp{B}(k)$.
	Otherwise, $\grp{G}$ has trivial intersection with the kernel of $\pi$,
	so it is explicitly isomorphic to 
	the subgroup $\pi(\grp{G})$ of $\grp{A}(k)$.
\end{proof}

\begin{corollary}
\label{Ext-hardest-DLP}
	The DLP in $\grp{C}(k)$
	is no harder than the hardest DLP in $\grp{A}(k)$ and $\grp{B}(k)$.
\end{corollary}

Proposition~\ref{DLP-reduction}
shows that if $\grp{G}$ is not contained in $\iota(\grp{B})$,
then the DLP in $\grp{G}$ reduces to the DLP in $\grp{A}(k)$.
It is important to note that the absence of a natural projection
from $\grp{C}$ to $\grp{B}$
does {\em not} preclude the existence of a homomorphism 
mapping $\grp{G}$ into $\grp{B}$;
thus the DLP in $\grp{G}$ may, in some cases,
be reduced to the DLP in $\grp{B}(k)$ as well.
For many subgroups $\grp{G}$, therefore,
the DLP in $\grp{G}$ is only as hard as 
the easier of the DLP in $\grp{A}(k)$ and the DLP in $\grp{B}(k)$.
This means that we can have a relative loss in security
in using extensions of $\grp{A}$ by $\grp{B}$
rather than using $\grp{A}$ and $\grp{B}$ independently.

\begin{remark}	
	Couveignes~\cite{Couveignes} shows that 
	if $\grp{C}$ is a commutative algebraic group extension of $\grp{A}$ by $\grp{B}$,
	then there exists an algorithm to solve the DLP in $\grp{C}$
	in subexponential time in the size of $\grp{C}$
	if and only if 
	there exists such algorithms 
	for $\grp{A}$ and for $\grp{B}$~\cite[Theorem 2]{Couveignes}.
	This is due to the existence of a $k$-rational isogeny 
	(not constructed in~\cite{Couveignes})
	from $\grp{C}$ to the direct product $\grp{A}\times \grp{B}$.
\end{remark}

\section{Extensions Presented by Cocycles}

Extensions $\grp{C}$ of $\grp{A}$ by $\grp{B}$ are effectively determined
by the choice of a symmetric $2$-cocycle
({\em cocycle} in the sequel):
that is, a map $c: \grp{A}^2 \to \grp{B}$
satisfying the relations
\begin{equation}
\label{cocycle-relations}
	c(P,Q) + c(P+Q,R) = c(Q,R) + c(P,Q+R) 
	\text{\quad and\quad }
	c(P,Q) = c(Q,P) 
\end{equation}
for all $P$, $Q$ and $R$ in $\grp{A}$.
Note that $c$ is {\em not} required to be a homomorphism.

Given a cocycle $c: \grp{A}^2 \to \grp{B}$,
we construct an extension $\grp{C}$ of $\grp{A}$ by $\grp{B}$ as follows.
The underlying variety of $\grp{C}$ is the direct product $A\times B$,
the identity element is $(0_\grp{A}, 0_\grp{B})$,
and the group law and inverse maps
are the morphisms $m_c: (A\times B)^2 \to A\times B$
and $i_c: A\times B\to A\times B$ defined by
$$
	m_c: 
	\left( (P_\grp{A},P_\grp{B}) , (Q_\grp{A}, Q_\grp{B}) \right) 
	\longmapsto
	\left( P_\grp{A} + Q_\grp{A}, P_\grp{B} + Q_\grp{B} + c(P_\grp{A},Q_\grp{A}) \right)
$$
and 
$$
	i_c:
	( P_\grp{A}, P_\grp{B} ) 
	\longmapsto 
	( -P_\grp{A}, -P_\grp{B} + c(P_\grp{A},- P_\grp{A}) )
$$
(here $+$ and $-$ denote group operations in $\grp{A}$ and $\grp{B}$).
Note that associativity and commutativity follow from
the relations~\eqref{cocycle-relations} above.
We say that $\grp{C}$ is the algebraic group 
{\em presented by the cocycle~$c$}.
Generalized Jacobians (for background, see \cite{Serre})
are examples of extensions presented by cocycles;
we will give an example below.
The direct product group $\grp{A}\times \grp{B}$
is the extension presented by the zero cocycle,
sending each element of $\grp{A}^2$ to $0_\grp{B}$.
Our assumption that $\iota$ and $\pi$ are explicit
holds in any extension presented by a cocycle,
as shown by the following easy lemma.

\begin{lemma}
\label{explicit-maps}
	Let 
	$
		0 
		\to 
		\grp{B} 
		\stackrel{\iota}{\to} 
		\grp{C} 
		\stackrel{\pi}{\to} 
		\grp{A} 
		\to 
		0
	$ 
	be an extension 
	presented by a cocycle $c: \grp{A}^2 \to \grp{B}$.
	
	\begin{enumerate}
	\item	The injection $\iota: \grp{B} \to \grp{C}$ 
		is given by $\iota(P) = (0_\grp{A}, P)$.
		The subgroups of $\grp{C}$
		in the image of $\iota$
		are precisely those of the form 
		$\{ (0_\grp{A}, P) : P \text{ in some subgroup of } \grp{B} \}$,
		and in such groups the map $\iota^{-1}$
		given by $\iota^{-1}((0_\grp{A},P)) = P$
		reduces the DLP to a DLP in $\grp{B}$.
	\item 	The projection $\pi: \grp{C} \to \grp{A}$ 
		is given by $\pi(P_\grp{A}, P_\grp{B}) = P_\grp{A}$.
		This map reduces the DLP in any subgroup of $\grp{C}$
		not in the image of $\iota$
		to a DLP in $\grp{A}$.
	\end{enumerate}
\end{lemma}

Lemma~\ref{explicit-maps} implies
that the difficulty of the DLP in $\grp{C}$
cannot be increased by a ``clever'' choice of $c$.
Indeed, each prime-order subgroup $\grp{G}$ of any extension $\grp{C}$
either projects faithfully into $\grp{A}$
or can be pulled back to $\grp{B}$.
In particular,
the DLP in any extension of $\grp{A}$ by $\grp{B}$
is no harder than the DLP in the direct product $\grp{A}\times\grp{B}$.

Suppose $\grp{C}$ is an extension 
presented by a cocycle $c: \grp{A}^2\to\grp{B}$.
Computing the group law in $\grp{C}$
requires the same computations
as computing the group law in $\grp{A}\times\grp{B}$,
together with an application of the cocycle $c$
and an extra group operation in $\grp{B}$
--- 
so computing the group law in $\grp{C}$
requires at least as much space and time as 
computing the group law in $\grp{A}\times\grp{B}$.
Further,
$\grp{C}$
and $\grp{A}\times\grp{B}$
have the same underlying variety,
so representing their elements requires the same space.
Thus computing in $\grp{C}$ requires at least as much time and space
as computing in $\grp{A}\times\grp{B}$.

For the purposes of DLP-based cryptography,
the group $\grp{A}\times \grp{B}$ 
offers no advantages over $\grp{A}$ and $\grp{B}$.
We have seen that
the DLP in $\grp{A}\times \grp{B}$ 
can be no harder than the hardest DLP in $\grp{A}$ or $\grp{B}$,
and computing in $\grp{A}\times\grp{B}$
requires at least as much space and time as 
computing in $\grp{A}$ and $\grp{B}$ separately.
Therefore,
using $\grp{A}\times \grp{B}$ in a DLP-based cryptosystem 
in place of $\grp{A}$ or $\grp{B}$
offers no advantage in security,
while requiring more storage space and computing time.
Similarly,
using an extension $\grp{C}$ presented by a cocycle $c: \grp{A}^2 \to \grp{B}$
instead of $\grp{A}$ or $\grp{B}$ alone
offers no increase in security,
since it has no larger prime-order subgroups
than those already present in $\grp{A}$ and $\grp{B}$,
while requiring at least as much time and space as 
computing in $\grp{A}$ and $\grp{B}$.
We have thus derived the following result.

\begin{proposition}
\label{summary-prop}
	If $\grp{C}$ is an extension of $\grp{A}$ by $\grp{B}$
	presented by a cocycle,
	then any DLP-based cryptosystem based on 
	a subgroup of $\grp{C}(k)$
	\begin{itemize}
		\item	is no more secure,
		\item	takes more space, and
		\item	is less computationally efficient
	\end{itemize}
	than the analogous cryptosystem based on $\grp{A}(k)$ or $\grp{B}(k)$
	(whichever has the harder DLP).
\end{proposition}

\begin{example}
\label{gec-example}
	In \cite{Dechene-thesis} and \cite{Dechene-ANTS},
	D\'ech\`ene proposes certain generalized Jacobians of elliptic curves
	as a supply of cryptographic groups.
	Suppose $\grp{E}$ is an elliptic curve over $k$,
	and let $O$ be the identity of $\grp{E}$.
	Let $\multiplicativegroup$
	denote the multiplicative group over $k$
	(we will write its group law multiplicatively).
	Fix points $M$ and $N$ (neither equal to $O$) on $\grp{E}$;
	the effective divisor
	$\modulus = (M) + (N)$
	is called the {\em modulus}.
	The generalized Jacobian $\Jac{\grp{E},\modulus}$
	is defined to be the extension of $\grp{E}$ by $\multiplicativegroup$
	presented by the cocycle
	$ c_\modulus(P,Q) = f_{P,Q}(M)/f_{P,Q}(N) $,
	where $f_{P,Q}$ is any function on $E$
	with divisor $(P + Q) + (O) - (P) - (Q)$.
	\footnote{
		We may take $f_{P,Q} = v/l$,
		where $l$ is the line through $P$ and $Q$,
		and $v$ is the vertical line through 
		the third point of intersection of $l$ with $E$.
	}
	The group law on $\Jac{\grp{E},\modulus}$ is given by
	$$
		(P,\lambda) + (Q,\mu) 
		= 
		(P + Q, \lambda\cdot\mu\cdot  c_\modulus(P,Q) ) .
	$$
We remark that this group law was also used in
Section 3 of \cite{Frey-Ruck}.
	
	The observations of Proposition~\ref{summary-prop} 
	all apply to $\Jac{\grp{E},\modulus}$.
	We know 
	that the DLP in $\Jac{\grp{E},\modulus}(k)$
	is no harder than 
	the hardest DLP in $\grp{E}(k)$ and $\multiplicativegroup(k)$.
	Thus using cyclic subgroups of $\Jac{\grp{E},\modulus}(k)$
	instead of subgroups of $\grp{E}(k)$ or $\multiplicativegroup(k)$
	requires extra work, and extra space, for no gain in security.
Indeed, it is widely recognised that elliptic curves give better
security and performance than multiplicative groups of finite
fields.
	Hence, it would be better either to remove the 
$\multiplicativegroup(k)$ and
	use only $\grp{E}(k)$ (saving space and time),
	or 
	spending the extra bits on a larger ground field $K$
	and using a prime order elliptic curve
$\grp{E}(K)$ instead (maximizing security).
\end{example}
	
\begin{remark}
	D\'ech\`ene suggests taking $M$ and $N$
	to be defined over a finite extension $K/k$,
	so that the cocycle $c_\modulus$ 
	maps $\grp{E}(k)^2$ into $\multiplicativegroup(K)$,
	and such that both $\grp{E}(k)$ and $\multiplicativegroup(K)$
	contains a subgroup of prime order $l$.
	Balasubramanian and Koblitz~\cite{Balasubramanian-Koblitz}
	have shown
	that for general elliptic curves,
	the degree of the smallest such extension 
	(called the {\em embedding degree})
	tends to grow with $l$,
	rendering computation in $\multiplicativegroup(K)$
	and $\grp{G}$
	exponentially difficult. 
	In practice,
	therefore,
	the suggestion requires
	$\grp{E}$ to be a so-called {\em pairing-friendly} curve,
        which means there is a homomorphism from 
        $\grp{E}$ to $\multiplicativegroup$ as used in
	the Frey--R\"uck and MOV attacks \cite{Frey-Ruck,MOV}.
	As a result, 
	this suggestion weakens $\grp{E}$,
	and therefore (by Corollary~\ref{Ext-hardest-DLP})
	weakens $\Jac{\grp{E},\modulus}$.
	In fact, as noted above, the generalised Jacobian
group law is the same as the method proposed by \cite{Frey-Ruck} 
for computing the Tate pairing (except the function is
inverted).
Hence, if $m$ is the least common multiple of the order
of $P$ and the order of $M-N$ in $\grp{E}(K)$,
then computing $m$ times $(P, 1)$ gives $( 0, \langle P, M-N
\rangle_m^{-1} )$.
\end{remark}

%

\end{document}